\documentclass[journal,twocolumn,10pt]{IEEEtran}

\pdfoutput=1

\usepackage{amsmath,bm}

\usepackage[]{times}
\usepackage{epsfig}
\usepackage{subfigure}
\usepackage{amsmath}
\usepackage{amssymb}
\usepackage{graphicx}
\usepackage{subfigure}
\usepackage{multirow}
\usepackage{flushend}
\usepackage{url}
\usepackage{color}
\newcommand{\bega}{\begin{eqnarray}}
\newcommand{\ega}{\end{eqnarray}}
\newcommand{\bb}{\begin{equation}}
\newcommand{\ee}{\end{equation}}

\begin{document}

\title{Predicting Failures in Power Grids:\\ The Case of Static Overloads}
\author{Michael Chertkov,
\thanks{M. Chertkov [chertkov@lanl.gov] is with Theory Division, CNLS \&
New Mexico Consortium,Los Alamos, NM 87545(4), USA.}
{Feng Pan},
\thanks{F.~Pan [fpan@lanl.gov] is with D-division, LANL, Los Alamos, NM 87545, USA.}
and Mikhail G. Stepanov
\thanks{M.~G.~Stepanov [stepanov@math.arizona.edu] is with Department of Mathematics, UA
Tucson, AZ 85721, USA.}}

\markboth{Submitted to IEEE Transactions of Smart Grids}{Chertkov, Pan, Stepanov: Predicting Failures in Power Grids}
\maketitle

\begin{abstract}
Here we develop an approach to predict power grid weak points, and specifically to efficiently identify the most probable failure modes in static load distribution for a given power network. This approach is applied to two examples: Guam's power system and also the IEEE RTS-96 system, both modeled within the static Direct Current power flow model. Our algorithm is a power network adaption of the worst configuration heuristics, originally developed to study low probability events in physics and failures in error-correction. One finding is that, if the normal operational mode of the grid is sufficiently healthy, the failure modes, also called instantons, are sufficiently sparse, i.e. the failures are caused by load fluctuations at only a few buses. The technique is useful for discovering weak links which are saturated at the instantons. It can also identify  generators working at the capacity and generators under capacity,  thus providing predictive capability for improving the reliability of any power network.
\end{abstract}

\begin{keywords}
Power Flow, Rare events, Distance to Failure
\end{keywords}

The basic features of the electrical power grid  have remained unchanged for over a hundred years. It is a hierarchical, centrally-controlled structure that assumes that power is generated solely from large central facilities, that power is abundant, and that power generation is relatively benign. Reliability is attained through redundancies and highly controllable generation that reacts to problems such as demand fluctuations and outages, rather than anticipating and avoiding them. A number of significant recent developments have upended these assumptions \cite{09EERE,09EDR}. Demand has risen sharply and become more variable in recent decades adding significant stress to the grid.  Alternative renewable energy sources, including wind and solar, are being integrated into the power grid but their distributed and intermittent nature make them difficult to control. With the existing grid control mechanisms and the ongoing developments, the probability of the existing power grid to survive failures or attacks degrades significantly,  emphasizing the need for new tools for power grid analysis.

Normally, the operation of the grid keeps demand and load in balance and generators under their capacities. A power-flow solver, simulating the grid according to the Kirchhoff laws and technological constraints \cite{96WW,07-Nag, 07-Glo}, finds acceptable generation and transmission capacity to meet consumer demand in this SATISFIABLE (SAT) state constrained by the physical properties of the grid. In the power engineering literature this SAT domain, including all the SAT states, is usually called the feasibility region. Starting from a SAT state and perturbing it, for example by increasing loads, one typically arrives at another SAT state.  Yet, there are marginal SAT states, defining an ERROR SURFACE in the space of consumer demand, where such infinitesimal perturbations can push the grid into an UNSATISFIABLE (UNSAT) state, where no feasible power-flow solution exists. Standard term used in the power engineering literature for `` the error surface" is
``the feasibility boundary". How do we determine which perturbations are the most dangerous and what regions of the grid must be shed for the grid to remain operational? This paper addresses this difficult challenge systematically.

We built our approach on a methodology originated in physical applications, specifically the physics of disordered media \cite{64Lif,66ZL}, field theory \cite{77Lip} and statistical hydrodynamics \cite{96FKLM}. The basic idea we borrow from these fields (which generally deal with a huge space of configurations each weighted according known probabilities) consists of relating the tail of a certain observable distribution function to a much smaller phase space,  or in the extreme case of weak fluctuations to special instance(s)/configuration(s) in the phase space - the instantons. Obviously,  the basic idea of dealing with a huge phase space of configurations (each characterized by an assigned weight, probability or fitness) is not unique to physics and can readily be extended to other disciplines. Moreover, the original theoretical approach suggests an algorithmic solution to the applied problem of ``finding a needle in a haystack",  i.e. not only estimating the tail of of a probability distribution function but more importantly, identifying those specific configurations that are part of the distribution function tail which cannot be sampled or analyzed efficiently via brute force Monte-Carlo techniques. The ``applied" power of the instanton approach was recently utilized in the field of error-correction, where it was used to guide efficient searches for rare erroneous configurations that degrade the performance of decoding algorithms \cite{05SCCV,08CS} that are successful under normal circumstances. Related algorithmic techniques were suggested for sampling rare trajectories in dynamical systems with nontrivial (and often topologically nontrivial) energy landscapes \cite{04EV,07TK}.

Even though the instanton approach we suggest here is new to the field of power system modeling,  the questions of accessing the distance to failure and probing rare but devastating events in power systems have been discussed many times in the past. The questions of the power grid reliability have received special attention recently as a result of large scale and cascading failures, such as the infamous 2003 East Coast blackout \cite{NERC}. Both microscopic \cite{05CTD} and macroscopic \cite{07DCLN} models assessing the probability of  blackout cascades were proposed and the results were compared with the electric power disturbance data. (See also \cite{CascadeReview} for a brief review of methods of cascading failure analysis.) It was also argued \cite{03TAILM} that prevention, based on an instantaneous probing of the most probable failure mode from the current state estimate, may be the only cure against the blackouts in the future. The research reported in this paper,  in spite of dealing with an instant failure and not a cascading failure, aims to achieve the larger goal of understanding and preventing cascading failures by identifying the ERROR SURFACE (feasibility boundary) between the normal and problematic regions of grid operation. In this regard, our approach complements other techniques, notably \cite{09SWB,08BV,10PMDL}, developed to access static vulnerabilities and the distance to failure in power grids.

The remainder of the manuscript is organized as follows.
In Section \ref{sec:DC} we introduce the Directed Current (DC) power flow model testing if the current state of the load is SAT or UNSAT.
Section \ref{sec:Distance} describes our probabilistic distance to failure approach. It contains three Subsections describing details of implementation
Section \ref{subsec:Implementation}, probabilistic interpretation
Section \ref{subsec:Prob},  and relation and differences between our scheme and other related techniques from power engineering
Section \ref{subsec:relations}. Section \ref{sec:sim} discusses numerical tests of the instanton scheme. Specifically we analyze two data sets of the island of Guam and of the IEEE RTS-96 system in Section
\ref{subsec:Guam} and Section
\ref{subsec:RTS-96} respectively.  We also explain interesting peculiarities observed in these realistic data sets on an illustrative example of a triangular graph discussed in Section \ref{subsec:triangle}. Section \ref{sec:conc} is reserved for discussions of Results, Conclusions and Path Forward.

\section{DC Power Flow }
\label{sec:DC}

We adapt a standard static model of the grid by considering it as a graph with nodes representing loads, generators, and control equipment and edges representing power lines/links. The basic equations governing the static realization of electrical power flow of electricity (a version of the Kirchhoff equations for currents and voltages) express the vector of power consumption and generation as a quadratic form of the vector of voltages, where coefficients in the quadratic form are complex impedances of different edges of the grid \cite{96WW,07-Nag, 07-Glo}. Solutions of the power flow equations should additionally satisfy the set of operational constraints describing currents, voltages, phases and other graph-specific constraints. The so-called DC (Direct Current) approximation
assumes that the resistivity-to-reactance ratio is small for any link, the voltage amplitudes are the same across the grid, and changes in the phase between neighboring nodes are sufficiently small.

Within this DC approximation we seek the optimal distribution of generation that keeps the amount of load shed over the given network to minimum. We state this optimization in terms of the following Linear Programming (LP), and computationally tractable, problem:
\begin{eqnarray}
 && LP_{DC}({\bm d}|{\cal G};{\bm x};{\bm u};{\bm P}) =\min_{{\bm f},{\bm \theta},{\bm p},{\bm s}}
 \left(\sum_{i\in{\cal G}_d} s_i\right)_
{COND},
\label{LP_DC}\\
 && COND = COND_{flow} \cup COND_{DC} \cup COND_{edge} \nonumber\\ && \cup
COND_{power} \cup COND_{over}, \label{COND}\\
 && COND_{flow}=\Biggl(\forall i\in{\cal G}_0:\quad \sum_{j\sim i}
f_{ij}\nonumber\\ && =\Biggl\{\begin{array}{cc}
 p_i, & i\in{\cal G}_p\\
 -d_i+s_i, & i\in{\cal G}_d\\
 0, & i\in{\cal G}_0\setminus({\cal G}_p\cup{\cal G}_d)\end{array}\Biggr.\Biggr),
 \label{flow_cond}\\
 && COND_{DC}=\Biggl(\forall \{i,j\}\in{\cal G}_1:\ \ \theta_i-\theta_j=x_{ij}f_{ij}\Biggr),
 \label{DC_cond}
 \end{eqnarray}
 \begin{eqnarray}
 && COND_{edge}=\Biggl(\forall \{i,j\}\in{\cal G}_1:\quad |f_{ij}|\leq u_{ij}\Biggr),
 \label{edge_cap_cond}\\
 && COND_{power}=\Biggl(\forall i\in{\cal G}_p:\quad 0\leq p_i\leq P_i\Biggr),\label{power_cap_cond}\\
 && COND_{over}=\Biggl(\forall i\in{\cal G}_p:\quad 0\leq s_i\leq d_i\Biggr),\label{over_cond}
\end{eqnarray}
where ${\cal G}_0$ and ${\cal G}_1$ are the set of vertexes and
edges, respectively, of the undirected power graph ${\cal G}$, and ${\cal
G}_p$, ${\cal G}_d$ stand for the subsets of generator and demand nodes,
i.e. ${\cal G}_p,{\cal G}_d\subseteq{\cal G}_0$. $i \sim
j$ means that the $i$ and $j$ nodes are connected by an edge, or $\{i, j\} \in
{\cal G}_1$. Here in Eqs.~(\ref{LP_DC}) we adopt standard notations for those vectors
associated with the set of vertexes and the set of directed edges, thus
 ${\bm f}=(f_{ij}=-f_{ji}|\{i,j\}\in
{\cal G}_1)$, ${\bm \theta}=(\theta_i|i\in {\cal G}_0)$,
${\bm p}=(p_i|i\in {\cal G}_p)$,
${\bm d}=(d_i|i\in {\cal G}_d)$, ${\bm s}=(s_i|i\in {\cal G}_d)$,
${\bm x}=(x_{ij}=x_{ji}|\{i,j\}\in {\cal G}_1)$,
${\bm u}=(u_{ij}\geq 0|\{i,j\}\in {\cal G}_1)$ and
${\bm P}=(P_{i}\geq 0|i\in {\cal G}_p)$
are vectors of power flows, phases, power generation, power demand, power shed,
reactance, flow capacity and maximum power generated, respectively.
As follows from Eq.~(\ref{COND}), the full set of
conditions is the union of sub-sets described in
Eqs.~(\ref{flow_cond}-\ref{over_cond}). These express respectively the
flow conditions, DC conditions on phases, capacity conditions (limits) on powers flowing along
the edges, generator capacity on the power the maximum amount of power the generators might be able to produce, and the natural conditions requiring that the positive amount of load
shedding at any consumer node cannot exceed the respective demand. Note
that inclusion of the DC conditions on the phases (\ref{DC_cond}) makes the
power grid flow optimization problem different from the standard LP formulation of the transportation flow optimization problem \cite{01CLRS}. For example controlling auto traffic on a graph of roads represents what $LP_{DC}$ would be without
Eq.~(\ref{DC_cond}). Our notation for the DC optimization scheme,
$LP_{DC}({\bm d}|{\cal G};{\bm x};{\bm u};{\bm p})$, emphasizes that in
the following we will consider variations of the load configuration,
${\bm d}$, while other characteristics, i.e. ${\cal G}$,
the vector of reactances, ${\bm x}$, the vector of power caps on edges,
${\bm u}$, and the vector of power generation capacity of the network,
${\bm P}$, remain fixed. Also, to simplify the notation, we use the following shorthand: $LP_{DC}({\bm d})$.

\section{Distance-to-Failure Algorithm}
\label{sec:Distance}

Optimization over distribution of generation, currents and  phases are done internally within $LP_{DC}({\bm d})$.
If a power flow solution which does
not require shedding exists, the $LP_{DC}$ outputs zero
in the normal SAT state. On the contrary, if a vector of
the demand, ${\bm d}$, outputs a nonzero  value (and thus positive by construction) for
$LP_{DC}({\bm d})$, we say that this ${\bm d}$ is UNSAT.

We introduce the distribution of the demand vector, ${\cal P}({\bm d})$, and aim to discover the so-called instantons, defined as configurations, ${\bm d}$, which are local maxima over the error-surface of  ${\cal P}({\bm d})$. We will be looking for multiple local maxima of the probability,  and especially for the absolute maximum, i.e. the most probable instanton, defined as the argument of $\left.\max_{\bm d} {\cal P}({\bm d})\right|_{ LP_{DC}({\bm d})>0}$. More generally we aim to output an (approximate) ordered list of the $k$-most probable instantons.

For concreteness, but without loss of generality,  we will consider solutions of $LP_{DC}$ over the following Gaussian distribution of demands \footnote{
Note that this particular form of the demand distribution of Eq.~(\ref{Gauss}), representing a crude model of reality (and as such ignoring effects of correlations between different consumers of electricity due to similar weather conditions, as well as variability in demand fluctuations representing non-uniformity of the demand forecast), plays an exemplary role in our distance to failure analysis. On the other hand,  we should warn the reader about possible sensitivity of some of the conclusions drawn from the instanton analysis illustrated below on the form of the demand distribution.},
\begin{equation}
{\cal P}({\bm d})
=\exp\left(-\frac{1}{2T}\sum_{i\in {\cal G}_d} (d_i/\bar{d}_i-1)^2\right)
\prod_{i\in {\cal G}_d}\left(2\pi T \bar{d}_i^2\right)^{-1/2} ,
\label{Gauss}
\end{equation}
where $\bar{d}_i$ is the average demand at load-node $i$, which we refer to as the
normal operating point of the grid.  The width of the demand
distribution at a node is conjectured to be proportional to the respective typical value, $\bar{\bm d}_i$. Starting from a normal operational point safely within
the SAT state, i.e. $LP_{DC}(\bar{\bm d})=0$, we study the sensitivity of the instanton set to moderate changes in $\bar{\bm d}$. In Eq.~(\ref{Gauss}), $T$ has the physical meaning of dimensionless dispersion parameterizing the inverse strength of the
relative fluctuations in demand, i.e. fluctuations in $d_i/\bar{d}_i$.

\begin{figure}
\centering
\includegraphics[scale=1.8]{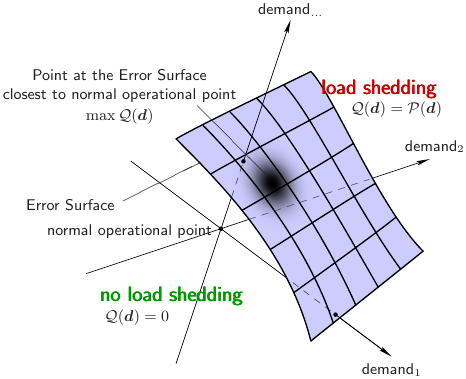}
\caption{\label{fig:smart_inst}The scheme illustrates a split of the space of demands into SAT (no shedding is needed) and UNSAT (shedding is unavoidable) domains. The error-surface is a surface in the space of demands separating the SAT and UNSAT domains. Normal operational point, ${\bm d}=\bar{\bm d}$, corresponds to the maximum of ${\cal P}({\bm d})$ defined in Eq.~(\ref{Gauss}). An instanton is a point at the error-surface where ${\cal P}({\bm d})$ achieves its maxima, or alternatively where $Q({\bm d})$, defined in the text, is maximized.  The most probable instanton, absolute maxima over possibly many instantons, is thus a point at the error-surface which is the ``closest" to the normal operational point.}
\end{figure}

The instanton approach was previously applied to error-correction in \cite{04CCSV,05SCCV}. The task of finding
the instantons is addressed by a similar method. We adapted this technique to the $LP_{DC}$
formulation. The instantons are formally defined as local maxima of the function
${\cal Q}(\bm{ d})$ such that ${\cal Q}(\bm{d}) = {\cal P}(\bm{
d})$, if the demands can not be met, and
${\cal Q}(\bm{d}) = 0$, if there is no shedding.
The function, ${\cal Q}(\bm{d})$, jumps from $0$ to a positive
value at the error surface, therefore guaranteeing that its maxima are achieved at the error-surface.  The decomposition of the phase space of demands into SAT/UNSAT domains and emergence of the error-surface and location of the most probable instanton are schematically illustrated in Fig.~\ref{fig:smart_inst}. Note, that the error-surface (feasibility boundary) has a complicated and generally non-convex shape,  thus resulting in multiple local maxima of ${\cal Q}(\bm{ d})$.

\subsection{Numerical Implementation}
\label{subsec:Implementation}

Following the scheme of
\cite{04CCSV,05SCCV}, we search for the maxima of ${\cal Q}({\bm
d})$ with the help of a general-purpose optimization technique,
specifically the downhill simplex (or ``amoeba'') method
\cite{AMOEBA,NumRec}. This method is well-suited for our setting because
the function, ${\cal Q}({\bm d})$, is not continuous
thus preventing us from using otherwise preferable gradient-type
optimization methods. The simplex method will allow us to adapt to the complex
landscape of the problem automatically without requiring any
additional input. (Notice,  that our choice of the implementation algorithm is subjective and was also guided by
consideration of simplicity. We envision using more sophisticated and ${\cal Q}({\bm d})$ specific algorithms in the future.)  The resulting optimization function is not
concave and we expect to find many instanton solutions. Each initialization of the
instanton-amoeba could lead to a new instanton. The initialization selects
an initial simplex (built on $N_d+1$ points in
the space of demands, where $N_d$ is the number of demand nodes) such that
${\cal Q}({\bm d})>0$ is nonzero for at least some points/vertexes of the simplex. In essence, the instanton-amoeba method evolves the simplex, via a sequential set of shifts, contractions and extensions, towards its
eventual collapse to a local maximum of ${\cal Q}({\bm d})$. Note that any modification of the simplex vertex requires evaluation of the $LP_{DC}$ at this new location in the demand space. Different random initiations will sample the space of instantons, thus generating the so-called instanton spectrum describing the frequency of a given instanton occurrence and also suggesting an estimation for the ordered list (with respect to their probability of occurrence and frequency) of top instantons.  Repeated an infinite number of times, the sampling would output the most probable instanton. However given that the number of initiations will be finite, the most probable one (of the finite number of
instantons found) gives an estimate from below for the probability of
the most probable instanton. To ensure sampling quality in our simulations
we continue random initiations till the most probable instantons would appear multiple number of times (typically, this required hundreds of initiations).

\subsection{Probabilistic Interpretation}
\label{subsec:Prob}

The instanton technique is solving inherently probabilistic problem in a deterministic way.  Indeed, when the relative fluctuations of demand are low, i.e. when the dimensionless dispersion, $T$, in Eq.~(\ref{Gauss}) is sufficiently small, the most probable instanton, ${\bm d}_{\mbox{inst};1}$, dominates the
asymptotic expression for the total probability of shedding (being UNSAT), ${\cal P}_{\mbox{shed}}\propto {\cal P}({\bm d}_{\mbox{inst};1})$. To describe the probability of shedding at  somewhat larger (but still not too large) values of the dimensionless dispersion parameter, $T$, one needs to account for other instantons,
\begin{eqnarray}
{\cal P}_{\mbox{shed}}\equiv \int{\cal Q}({\bm d}) {\cal D}{\bm d}
\approx \sum_{k}\mu_k{\cal P}({\bm d}_{\mbox{inst};k}),
\label{many_inst}
\end{eqnarray}
where ${\cal D}{\bm d}$ stands for the differential over $N$-dimensional space of demands and ${\cal P}_{\mbox{shed}}$ represents the integral of the probability density of ${\bm d}$ over the UNSAT part of the space of demands; $\mu_k$ accounts for the ``volume" of the respective instanton, i.e. for effect of UNSAT configurations which are close to ${\bm d}_{\mbox{inst};k}$.
It is also assumed in Eq.~(\ref{many_inst}) that different terms on the right hand side are ordered in $k$ according to their probability of occurrence, i.e. for any $k_2>k_1\geq 1:\quad{\cal P}({\bm d}_{\mbox{inst};k_2})<{\cal P}({\bm d}_{\mbox{inst};k_1})$. Thus, the volume of the shadowed domain in Fig.~(\ref{fig:smart_inst}) corresponds to
$\mu_1$ for the most probable instanton. One expects that, of the two factors under the sum in Eq.~(\ref{many_inst}), the dependence of $\mu_k$ on $1/T$ is much weaker (algebraic) than the fast (exponential) decay of ${\cal P}({\bm d}_{\mbox{inst};k})$ with increasing $1/T$.
Calculating $\mu_k$ would require an accurate exploration of the $k$-instanton neighborhood according to ${\bm Q}({\bm d})$, understood as an un-normalized sampling probability. We do not do these more involved computations in this paper.  However, we suggest that the frequency of a given instanton appearance within the instanton-amoeba procedure can be considered as a rough proxy for $\mu_k$.

\subsection{Relations and Differences with Other Approaches in Power Engineering}
\label{subsec:relations}

Let us briefly discuss some relations, but also differences, between our instanton approach and two important concepts in power engineering: first the Available Transfer Capability (ATC), see e.g. \cite{08WSI,97Sau,06Sch}, and second the probabilistic and stochastic Optimal Power Flows (OPF) \cite{06Sch}. ATC is normally defined as an available power transport reserve between two areas/regions connected by a link (possibly aggregated link). To compute ATC one increases flow through the link incrementally, introduce one (or few) degrees of freedom to adjust (e.g. overall rescaling of the generation within a slack bus model of generation) and test if the resulting power flow solver outputs solution free of capacity violations over the entire network considered. (In a more sophisticated version evaluation of ATC may also account for $N-1$ contingencies,  thus requiring that the resulting solutions do not violate capacity constraints even when any single element of the system is removed.) As such ATC also measures distance to failure,  however in a different sense then the instanton technique.  First of all, our instanton technique explores multi-dimensional variations (in the space of demands) while ATC normally tests a single parametric changes (in power flow along a single link). Also, and more importantly, the instanton technique is inherently probabilistic,  moreover focusing on predicting rare configurations resulting in a failure (these taken place with probability significantly less than unity). On the contrary,  even when ATC is discussed in a stochastic setting (for example assuming averaging over fluctuations of loads) the fluctuations considered are typical,  i.e. these realized with a reasonable, $O(1)$, probability. Similarly,  the concepts of probabilistic and stochastic OPF,  defined as a minimum cost generation optimized over spatial distribution of available generation and averaged over fluctuations (for example over probability distribution of demands), accounts only for typical fluctuations. The OPF is also similar to the $LP_{DC}$ scheme of Eq.~(\ref{LP_DC}) in what concerns optimization over spatial power distribution,  however in the probabilistic instanton search scheme the $LP_{DC}$ is an important but only intermediate step of evaluation.

\section{Simulation Results}
\label{sec:sim}

\subsection{Example of the Power Grid of Guam}
\label{subsec:Guam}

\begin{figure}[h]
\centering
\psfig{file=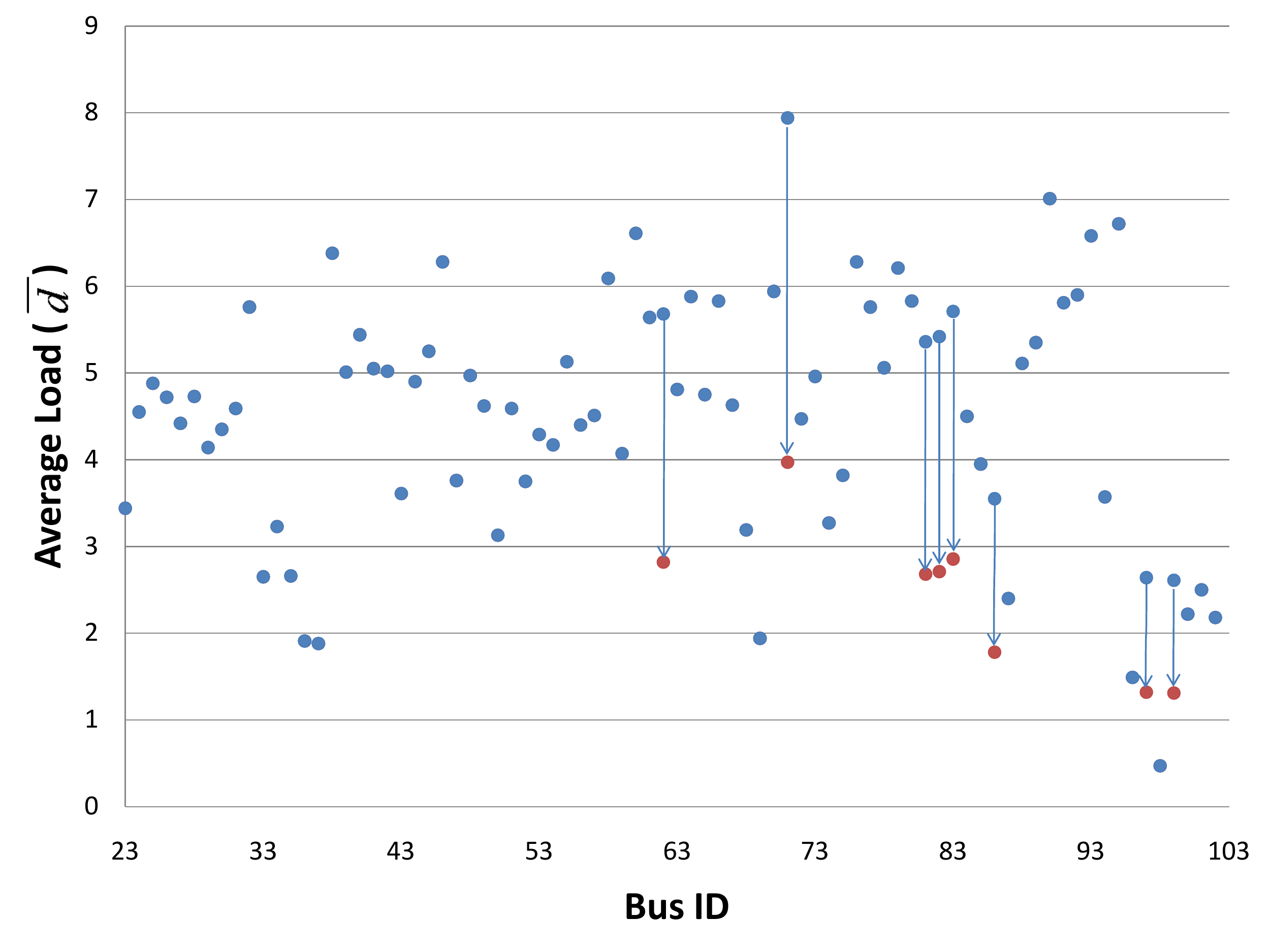, page=5,scale=.35}\\
\caption{\label{fig:Guam} \small
This Figure shows the power grid/graph of the island of Guam. Green and Yellow/Red/Orange/LightBlue nodes represent generators and loads, respectively.
Edges represent power lines. Parts of the graph (vertexes and edges) marked in Red/Orange/LightBlue/Black show localization/saturation of the three most probable instantons found via the instanton-amoeba algorithm described in the text. Edges shown in black are saturated on all the three instantons. Note,  that positions of the nodes and edges on the graph do not correspond to their actual geographical location (See the geographical map of Guam in the offset on the top of the Figure).
}
\vspace{-0.2cm}
\end{figure}

\begin{figure}[h]
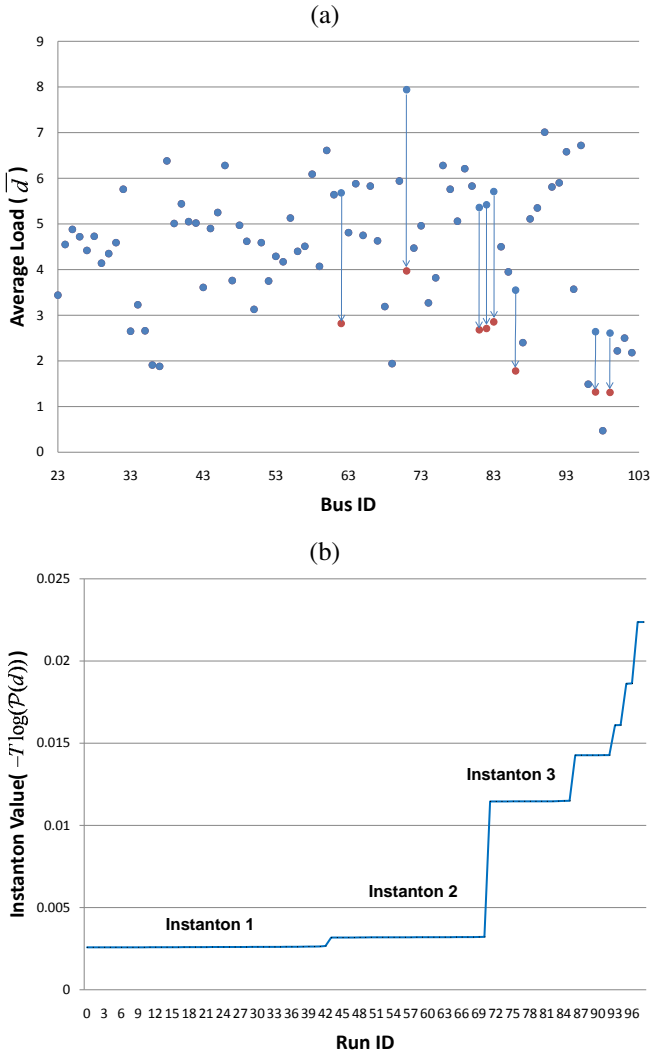

\centering
(a) \psfig{file=allfigures090710.pdf, page=1,scale=.35}\\
(b) \psfig{file=allfigures090710.pdf,page=2,scale=.35}
\caption{\label{fig:guam_demand1} \small
Demands and instanton values in the Guam grid.
Fig.~(a) shows the average load, $\bm d$, of the Guam grid (blue dots). The reduced average load configuration (see the main text for details) coincides with the original configuration
at all but the eight load nodes shown in red.
Fig.~(b) shows the demand distribution, $-T\log({\cal P}({\bm d}))$, for $100$ instanton-amoeba runs resulting in $8$ distinct instantons (only $7$ of which are shown in (b)).}
\vspace{-0.2cm}
\end{figure}

\begin{figure}[h]
\centering
\psfig{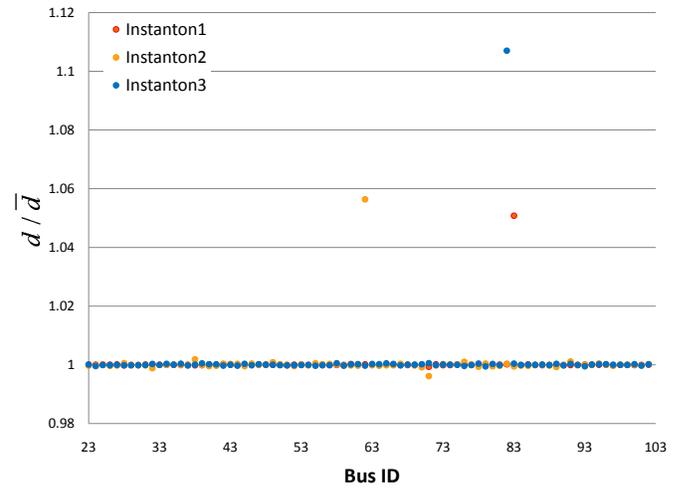}
\caption{\label{fig:guam_loads1} \small
Relative demand, $d_i/\bar{d}_i$, vs. $i$ in the Guam grid is shown in different colors for the $3$ most probable instantons.}
\vspace{-0.2cm}
\end{figure}

\begin{figure}[h]
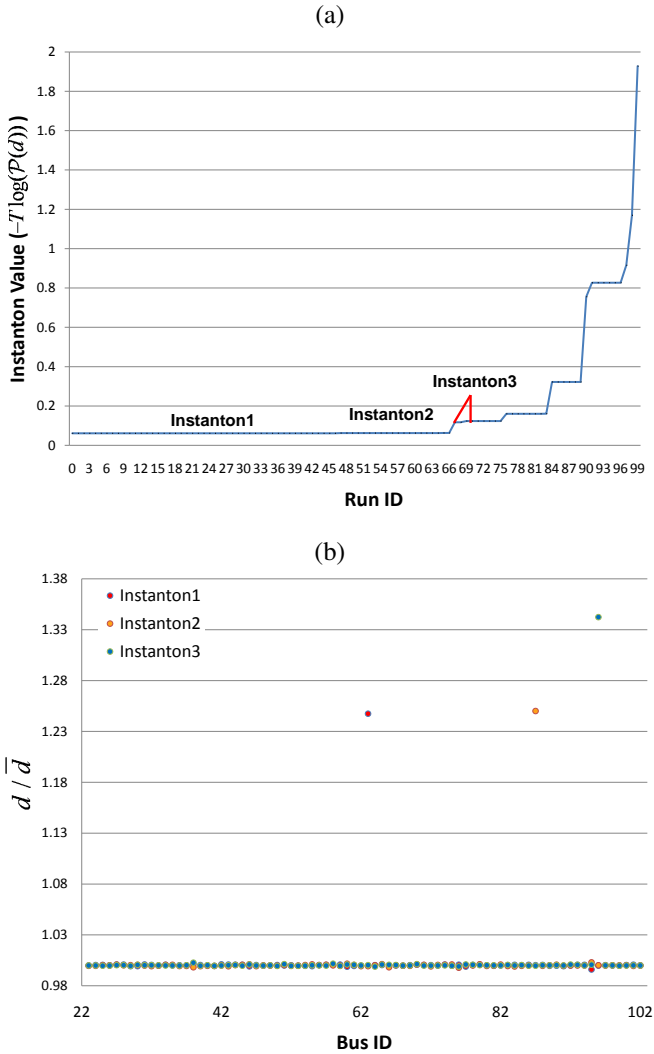

\centering
(a) \psfig{file=allfigures090710.pdf,page=6,scale=.35}\\
(b) \psfig{file=allfigures090710.pdf,page=8,scale=.35}
\caption{\label{fig:guam_reduced} \small
(a) Shows probability distribution of instantons, $-T\log({\cal P}({\bm d}))$, in the Guam example for the reduced average load (reduced values are shown in red in Fig.~\ref{fig:guam_demand1}a). One observes $10$ distinct instantons in $100$ runs of the instanton-amoeba. The instanton values of the top two instantons (observed
in $66$ total runs) are very close. Three instances, corresponding to third-ranked instanton,
are marked by red lines limiting the range. (b) This figure shows the three top-ranked instantons (in different colors) super-imposed in one plot.
}
\vspace{-0.2cm}
\end{figure}

In the first example, the instanton search algorithm is tested on the power grid system of the island of Guam, with the data received from a LANL infrastructure database. As seen in Fig.~\ref{fig:Guam}, The power grid of Guam consists of $103$ buses  ($23$ generators and $80$ loads) and $116$ lines/edges.  For average load, $\bar{\bm d}$, in distribution (\ref{Gauss}),  we use the actual values for a typical day illustrated in Fig.~\ref{fig:guam_demand1}(a).  The Guam grid model specifies  the line capacity vector, ${\bm u}$, and reactance, $\bm x$.
There were no generation capacity data available. We set the generation capacity to infinity, effectively eliminating condition (\ref{power_cap_cond}).

Our initial instanton-amoeba test for this system included $100$ runs. Fig.~\ref{fig:guam_demand1}(b) shows the distribution of  the instanton value, $-T\log(P({\bm d}))$, i.e. the negative logarithm of the instanton probability of occurrence measured in the units of fluctuation in demand intensity, $T$ \footnote{We included only $99$ values in Fig.~\ref{fig:guam_demand1}(b), as the $100$-th was seriously less probable and thus off-scale.}.
 The subset of the top-three instantons ordered with respect to their probability of occurrence is shown in Fig.~\ref{fig:guam_loads1} for comparison. One observes that the relative demand at most buses within each instanton is very close to the average demand (the normal operational point). For most frequently occurring instantons from the ordered list, there is a significant deviation only at one node or sometimes (for less probable instantons) two nodes. In other words, the instantons are remarkably sparse. Overall, there are $8$ problematic buses at which the resulting values of the demand, ${\bm d}_{\mbox{inst};k}$, are different from their respective values at the normal operational point. Three of the eight buses appear as anomalous in  Fig.~\ref{fig:guam_loads1}, and these three buses contributed significant deviations from the mean in $86$ of the $100$ instantons.
One of the most problematic buses (in the sense of its deviation value from the normal operational point in the multiple instantons)  is  $\# 83$, shown in Fig~\ref{fig:guam_loads1}. It contributed to $43$ (of the total of $100$) runs of the instanton-amoeba method.
Buses where the top three instantons from Fig.~\ref{fig:guam_loads1} show significant deviation from the mean value are shown in color in Fig.~\ref{fig:Guam}. We also show in Fig.~\ref{fig:Guam} edges for which the respective (top three) instanton solutions are saturated, in the sense that the output of the $LP_DC$ is realized with $COND_{edge}$ from Eq.~(\ref{edge_cap_cond}) being equality. Note that the top three instantons share among themselves $18$ common saturated lines (shown in black in Fig.~\ref{fig:Guam}), which are already saturated at the normal operational point,  i.e. at ${\bm d}=\bar{\bm d}$. Each of the instantons also has one additional new saturated line (shown in the proper color in Fig.~\ref{fig:Guam}) which appears attached to the problematic node  where the instanton is localized.

As seen in Fig.~\ref{fig:guam_loads1}, the instanton $\#2$ of the Guam model shows a counterintuitive feature: some loads need to be reduced, although slightly, in comparison with the normal operational point. To elucidate and explain this interesting feature related to correlations, imposed on loopy power grids via a combination of the phase constraints (\ref{DC_cond}) and link capacity constraints (\ref{edge_cap_cond}), we discuss below  how a similar behavior appears in a three-node (triangle) example.

Localization and sparsity of the observed instantons translates into  discovery of the weak point of the grid with respect to the given distribution of loads. It is suggestive to consider shifting the distribution of loads along the special dangerous directions in order to relieve the stress imposed on the system. (We assume here that redistribution of loads within the grid is possible.) To elucidate this point, we have reduced the average
demand at $5$ out of $8$ problematic buses by half, thus modifying the normal operational point of Fig.~\ref{fig:guam_demand1}(a). (The five positions and respective reductions are
shown in red in Fig.~\ref{fig:guam_demand1}(a). The choice of the number of nodes to be modified,  five in the example, is arbitrary and meant to illustrate the tendency.)
The instanton-amoeba search procedure for this new (shifted with respect to the
original operational point) distribution produced a set of new instantons, shown in Fig.~\ref{fig:guam_reduced}.
The top instantons for Fig.~\ref{fig:guam_reduced} are less probable than the top instantons in the original test, thus making the system less stressed. The configurations of the top three most probable instantons for the reduced value of demand are shown in Fig.~\ref{fig:guam_reduced}(b). Each of the top three instantons is localized on a bus, however after the targeted load reducing we find a new set of troublesome buses  distinct from those identified in the original demand configuration. The new buses can be viewed as the second least reliable set requiring attention.  We note that some of the higher (i.e. less probable instantons) become much less sparse/localized, although this is not demonstrated in the Figures. This iterative process of changing the mean load in the distribution (\ref{Gauss}) can obviously be continued to discover an ordered set of problematic buses.

\subsection{Example of the IEEE standard RTS-96 system}
\label{subsec:RTS-96}

\begin{figure}
\centering
\psfig{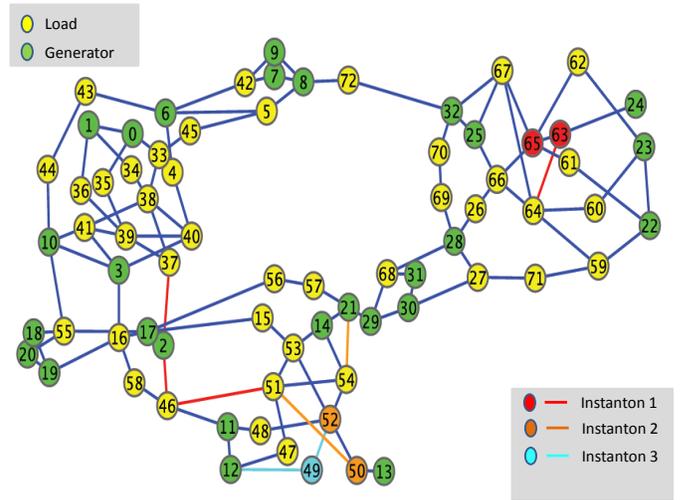}
\caption{\label{fig:RTS-96} \small
This Figure shows the power grid/graph of the IEEE RTS-96 system \cite{79-RTS,96-RTS} with color-coding and notations identical to these used in Fig.~\ref{fig:Guam}.
}
\vspace{-0.2cm}
\end{figure}

\begin{figure}
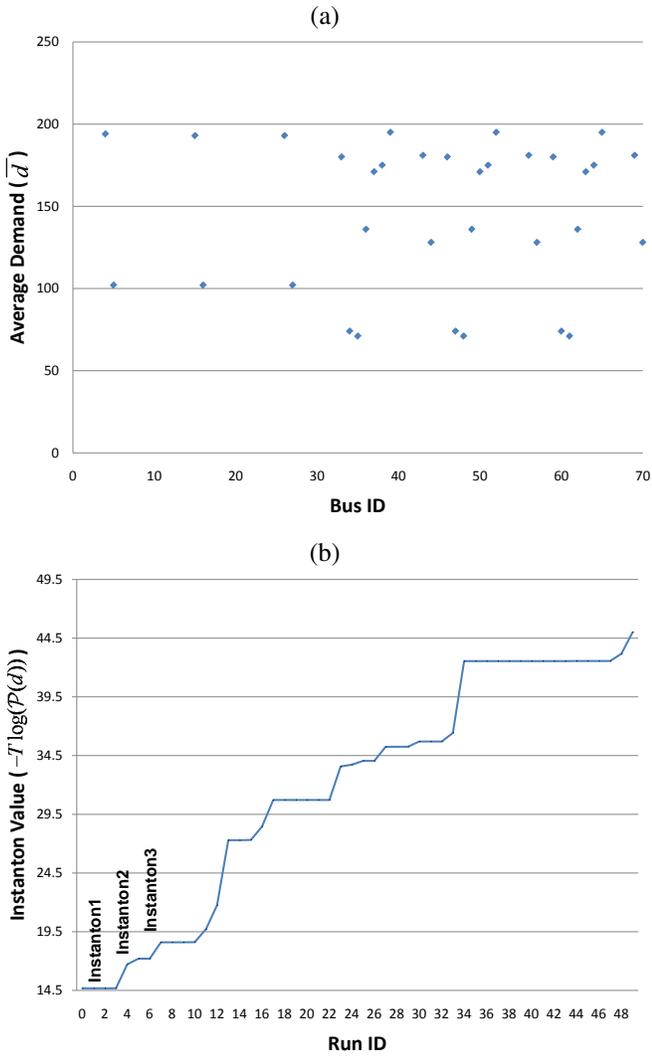

\centering
(a) \psfig{file=allfigures090710.pdf, page=10,scale=.35}\\
(b) \psfig{file=allfigures090710.pdf,page=11,scale=.35}
\caption{\label{fig:RTS_demand} \small
Normal operational point (average demand) and instanton distribution for example of the IEEE RTS-96 grid.
(a) shows our choice of the average $\bar{\bm d}$, (b) shows the load distribution ($-T\log({\cal P}({\bm d}))$ vs instanton ID) for $50$ runs of the instanton-amoeba resulting in $23$ distinct instantons.}
\vspace{-0.2cm}
\end{figure}

\begin{figure}
\centering
\psfig{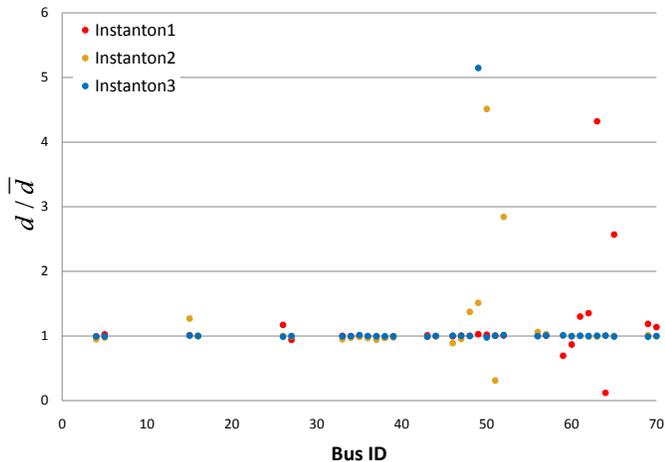}
\caption{\label{fig:RTS_loads} \small
Relative demand $d_i/\bar{d}_i$ vs $i$ in the RTS-96 grid for the top $3$ instantons.
}
\vspace{-0.2cm}
\end{figure}

\begin{figure}
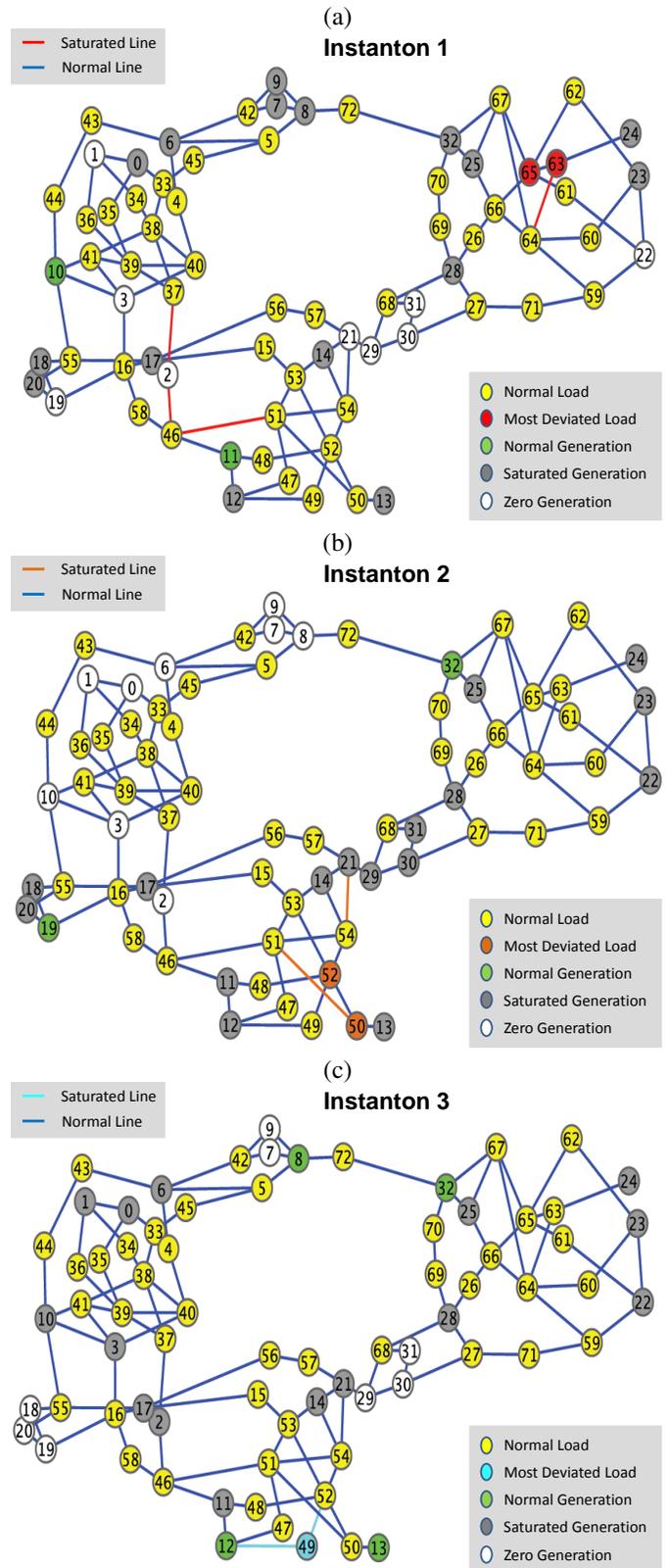

\centering
(a) \psfig{file=allfigures090710.pdf, page=15,scale=.35}\\
(b) \psfig{file=allfigures090710.pdf,page=16,scale=.35}
(c) \psfig{file=allfigures090710.pdf,page=17,scale=.35}
\caption{\label{fig:RTS_gen} \small
These Figures show feasible distributions of generation for the three leading instantons in the RTS-96 grid. The color coding for loads and lines is equivalent to the one used in Fig.~\ref{fig:RTS-96} and Fig.~\ref{fig:RTS_loads}. The generator nodes are colored here in green, gray and white depending on whether or not the generation was normal (within the bounds), saturated (at the upper bound), or was not utilized at all (zero generation) at the respective instanton configuration. Note,  that distribution of generation shown on the figure is not unique for the given instanton configuration.}
\vspace{-0.2cm}
\end{figure}

We also tested our instanton search algorithm on the IEEE-RTS-96 system and data set. The
detailed specifications of these data are described in \cite{79-RTS,96-RTS}. As seen in Fig.~\ref{fig:RTS-96}, the test system consists of $73$ buses and $120$ lines.  The generation capacity is distributed among $27$ buses, $33$ of the busses are dedicated to loads and the remaining $13$ are transmission buses.
The max load data is from \cite{96-RTS}  and we set the average load (also the main operational point) $\bar{\bm d}$ to half of the maximum load. Fig.~\ref{fig:RTS_demand}(a) shows the demand of the main operational point
\footnote{We originally ran the simulation choosing the main reference point to be at the maximum load point. (This maximum load point, correspondent to a very stressed day, was the only one available in the data set of \cite{79-RTS,96-RTS}). However, the multiple trials of the instanton-amoeba method then resulted almost always in one special instanton configuration. We have conjectured that this degeneracy was due to the fact that the maximum load point was very close to the error-surface, thus making the landscape searched by amoeba less rugged than in a more reasonable (less stressed) reference point resulting in a variety of competing instantons/minima.  Shifting the average load point to the half of the maximum load  indeed removed the redundancy. Note also, that we have modified the ordering of busses in comparison with the original ordering of \cite{96-RTS} to provide a better presentation in the charts.}.

Data from $50$ runs of the instanton-amoeba algorithm are shown in Fig.~\ref{fig:RTS_demand}(b). In this case the normal operational point is sufficiently far from the error-surface which shows itself in much smaller values of the instanton probabilities. We also observe more diversity. There are $22$ distinct instantons,  observed in $50$ runs, as opposed to $8$ distinct instantons observed in $100$ runs of our original Guam test. In contrast to the Guam test,  the instantons  in Fig.~\ref{fig:RTS_loads} are not as sparse.
However in Fig.~\ref{fig:RTS_loads} we still observe a significant localization phenomenon among the most probable instantons shown, although less pronounced. Instanton \#1 (observed $4$ times) is localized at buses $65$ and $63$,  where the instanton demand deviates significantly from the respective average load. Fig.~\ref{fig:RTS-96} also shows saturated edges for the top three instantons,  i.e. the edges where
the conditions Eqs.~(\ref{edge_cap_cond}) are tight (equalities and not inequalities).
Not all saturated edges in this example are graph-neighbors of the localized nodes, thus emphasizing the graph-global nature of the instanton search algorithm. This lack of graph proximity between the saturated edges and the localized nodes is also seen in instanton \#2 (observed only once in $50$ trials), but to a lesser degree. Generation data corresponding to the top three instantons are shown in Fig.~\ref{fig:RTS_gen}. ( Note that the generation data are not unique and allow slight freedom in redistribution.) We observed that the state of generation in these critical instanton configurations of the demand is really pushed to the limit and stressed. Indeed,  there are relatively few generators which operate normally (neither overstressed nor understressed) and the majority of generators are at capacity (shown in gray in Fig.~\ref{fig:RTS_gen}) or are not utilized at all (shown in white in Fig.~\ref{fig:RTS_gen})
\footnote{Note, that in setting the lower bound for all the generators to zero, we are ignoring important practical limitations associated with the fact that generation which is too small is not economically viable. We acknowledge that accounting for this important effect may have a significant impact on the distance to failure analysis, as
for example was proven to be important in the interdiction analysis of \cite{08BV}.}. We also observe that the generators performing at the capacity and significantly under capacity for the instantons are often grouped in graph-local cliques.  These cliques vary significantly from instanton to instanton, and the under capacity cliques typically lie on a path from the domains which are at the capacity (generation wise) to the most deviated loads \footnote{Note that such a large variability in generation is a consequence of the lack of a generator related contribution to the cost function of the optimization problem considered, $LP_{DC}({\bm d})$. A more accurate and realistic  modeling of the (operator controlled) cost function should be able to reduce this undesirable generation variability both in typical and atypical (instanton) configurations.}.

The feature of load reduction at some nodes (in comparison with values at the normal operational point), already noticed in instanton $\#2$ of the basic Guam model from Appendix \ref{subsec:Guam}, is also observed in the RTS-96 model, as seen in Fig.~(\ref{fig:RTS_loads}) in instantons $\#1$ and $\#2$. We explain this feature of meshed (i.e. containing many loops) power grids on a simple triangle example in the next Subsection.

Repeating the procedure of reducing the normal operational point at the localized nodes,  we have generally confirmed reduction of the instanton probability and emergence of the ``next most critical" elements of the grid, similar to the example of Guam. However,  we have also found a new and interesting situation in which we can force some instantons into SAT configurations by increasing demands at some nodes with instanton values of demand below respective mean values.
We suggest that this phenomenon is related to the global nature of correlations in power grids, formally contained in the conditions on phases (\ref{DC_cond}) of $LP_{DC}$. Indeed, these phase conditions in combination with edge-capacity constraints (\ref{edge_cap_cond}) may lead to selection of a relatively restrictive sub-space of possible transportation flow solutions described in Eq.~(\ref{flow_cond}). The restriction resulting the transition from SAT to UNSAT with decrease in demand, can be illustrated on a simple three-node example discussed below, and in fact, is the major factor behind the proposal of \cite{09HOFO}, suggesting that a smart use of switching can significantly extend the SAT domain by effectively removing some of the ``most restricting" phase conditions (\ref{DC_cond}).

\subsection{Triangle example}
\label{subsec:triangle}

\begin{figure}
\centering
(a) \psfig{file=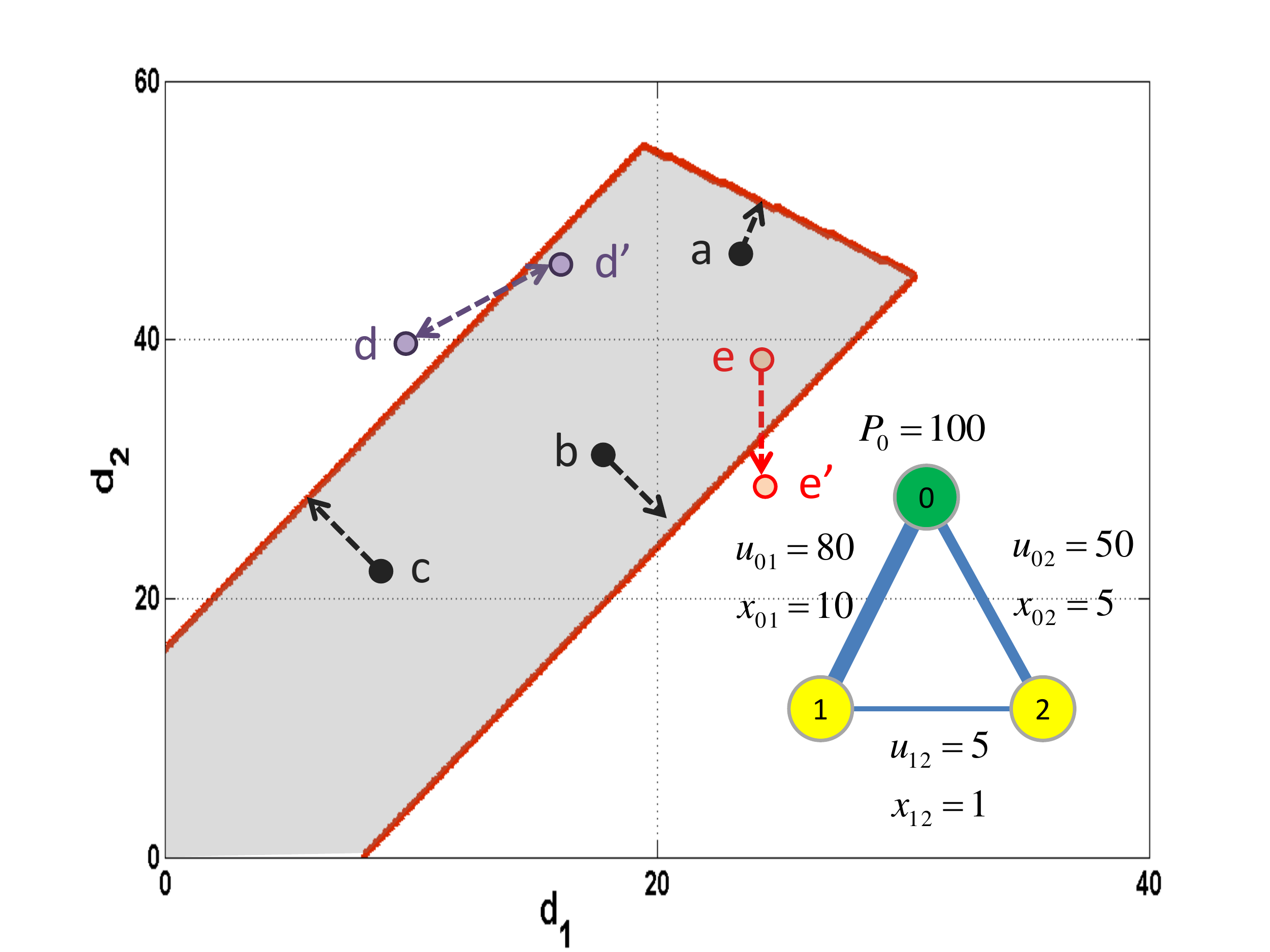,page=1,scale=.35}\\
(b) \psfig{file=triangle.pdf,page=2,scale=.35}
\caption{\label{fig:triangle}
(a) $(d_1,d_2)$ plane of demands of the triangle grid model shown in the lower right corner. Line and generator capacities, as well as reactances of the lines are fixed according to the numbers shown.
The gray area with the red envelope corresponds to the SAT domain and the area on the exterior of the red envelope marks the UNSAT domain. (a-f) points in the two-dimensional space of demand correspond to different representative configurations discussed in the text.
(b) More details on the distribution of currents, generation and load/shedding for (a-e) states of the triangle model shown in Fig. (a). Arrows and numbers next to arrows correspond to directions and values of currents, $f_{ij}$, over respective links. Numbers in boxes show generation/load, $p_i/(d_i-s_i)$, where the correction shown in red indicates the amount shed, $s_i$ ((d) is the only UNSAT state where shedding applies).}
\end{figure}

As observed above for the Guam and RTS-96 models, instantons may show reduction in some loads in comparison with respective normal operational points. Even more surprisingly, increase of some stressed loads in the basic operational point in the RTS-96 model (not modifying other loads,  thus corresponding to an overall increase in load) leads to a SAT state. Those features of the meshed power grid can in fact be understood on the basis of a simple triangular grid.

Consider a triangle with one generator, marked as ``0", connected to two loads,  marked ``1" and ``2". Assume that the $(1,2)$ link is significantly weaker than the capacities of the generator-to-load links and the generator capacity exceeds $d_1+d_2$. To elucidate the important consequence of this asymmetry in link capacities, let us also assume that characteristics of the $(0,1)$ and $(0,2)$ links are identical and the demands at the nodes $1,2$ are equal. (These restrictions will be removed later.) This state of the system demands is obviously SAT and it also corresponds to zero flow along the $(1,2)$ line between the two loads, with thus equal phases realized at the nodes $1$ and $2$. Now consider breaking the symmetry between the two loads by decreasing the demand at one of the loads, say by $1/2$, i.e. moving to the case with  $d_2=d_1/2$. Obviously this case becomes UNSAT if the capacity of the $(1,2)$ link, $u_{12}$, is sufficiently small.

This transition from SAT phase to UNSAT phase with a decrease in load demand is illustrated in Fig.~\ref{fig:triangle}a for a general triangle example (not requiring the generator-to-load links to be identical).
Details of the model (the graph, the capacities and the reactances) are shown in the lower right corner of Fig.~\ref{fig:triangle}a. The area shown gray/white in the two-dimensional demand space $(d_1,d_2)$ is SAT/UNSAT.
The SAT domain is convex. Different points,  marked (a-e), inside and outside of the domain illustrates counter-intuitive effects
encountered above in the Guam and RTS-96 models. Indeed, if point (b) or (c) is chosen for the main operational point, $\bar{\bm d}=(\bar{d}_1,\bar{d}_2)$, and the distribution of the demand is $-T\log({\cal P}(d_1,d_2))=(d_1-\bar{d}_1)^2+(d_2-\bar{d}_2)^2$, then transition to the most probable instanton, shown in Fig.~\ref{fig:triangle}a with arrow, corresponds to a decrease in one of the loads, $d_2$ or $d_1$. The normal operational point positioned at the point (a) has the most probable instanton with both loads increased. Moreover, selection of an UNSAT state (d) and simultaneous increase of both loads, as shown by respective arrows in Fig.~\ref{fig:triangle}a, results in a SAT state (d'). A simultaneous decrease of both loads transforms the SAT state (d') to UNSAT state (d). Finally, the decrease of
the demand $d_2$ while keeping $d_1$ constant transforms the SAT state (e) to UNSAT state (e'). Fig.~\ref{fig:triangle}a provides details on the distribution of currents and generation/load in the (a-e) states of the triangle model.

\section{Results, Conclusions and Path Forward}
\label{sec:conc}

Summarizing, the main findings of this work are:
\begin{itemize}
\item For a given power grid network considered within the DC power flow
approximation and under fixed configuration of demand, generation
capacity and transmission constraints, we have formulated the problem of
minimizing the load shedding over the network.

\item By analogy with coding theory,  we have introduced the error-surface, which is a surface in the space of loads that separates the region in which no load-shedding is required, called SAT, from the region in which some load-shedding is unavoidable, UNSAT.

\item Assuming load statistics with the typical load value being safely within the SAT region, we have posed the question of finding the most probable configuration in the multi-dimensional demand space which requires load shedding. More generally, we define the notion of instantons -- configurations on the error-surface correspondent to local maxima of the load probability. The problem of finding the instantons was stated as a generally non-concave maximization problem, in which multiple instantons are the respective local maxima.
\end{itemize}

This amoeba-based instanton-search scheme for the distance to failure was tested on two
exemplary networks shown in Fig.~\ref{fig:Guam} and Fig.~\ref{fig:RTS-96} respectively.
Here, we summarize interesting features of our analysis:
\begin{itemize}

\item If the typical value of the demand is sufficiently far from the error-surface, the most probable instantons are well-localized on the grid, in the sense that only few components of the resulting vector of demand are comparable to the respective values at the normal operational point.  The statement is even stronger for the ``less meshy" Guam network, in which the resulting relative deviation vectors for instantons are sparse.

\item We rank the instantons according to their probability of occurrence and find that the most probable instantons collectively are localized on a small number of nodes. In other words, there are normally a handful of problematic nodes that need to be watched for load shedding. Similarly,  we found that the number of problematic links that are saturated at the instanton solutions is relatively small. These links, as well as generators performing at the capacity, are actual reasons for the instantons/contingencies.  If capacities of the links and the generators are increased the respective instanton is removed from the list thus leading to the overall reduction in the failure probability of the system.

\item Introducing an upper bound on generation in one of the samples studied (RTS-96),  we observed that the number of generators at the capacity  and these significantly under capacity (i.e. generating negligibly small in comparison with the capacity) at the instanton configurations is significant, indicating the extremality of the stress experienced by the system.

\item We observed that the normal variation in $\bar{\bm d}$ does not drastically change the structure of the most probable instantons.  However, occasionally even a minor shift of $\bar{\bm d}$ may lead to qualitative changes in the structure of the error-surface (and thus change the resulting instantons). In particular, we observed that a reduction in some loads at the instanton state, in comparison with their values at the normal operational point state, is possible. This phenomenon is related to nonlocal correlations characteristic of power flows. Another phenomenon  of related origin, that we observed and explained on a simple triangle-grid example, consists of  possible transitions from a SAT state to an UNSAT state, under a decrease of one, many or possibly even all loads.

\end{itemize}

We view this work as the beginning of exciting further exploration of the power network stability along the following lines.
\begin{itemize}
\item In this study we have used the $LP_{DC}$ formulation as the simplest SAT-UNSAT criterion expressing specificity of the power grid. However any other and richer criteria (for example one representing an AC flow optimization or another one checking for dynamic instability of the grid)  can be used as a substitute for $LP_{DC}$ in  the elementary (inner-loop) step of the instanton scheme. In this regard, the instanton optimization scheme constitutes a black box capable of satisfying any sensible criterion.

\item On the other hand, the universality of the instanton-amoeba approach (as the outerloop of the scheme) can also be viewed as a handicap when it comes to discussing its complexity and implementation speed.  Thus,  in the field of error-correction, a significantly faster instanton-search scheme, customized to the specific form of the inner-loop criterion (LP-decoding, also the specific type of the channel statistics), was suggested in \cite{08CS}. This resulted in an essential acceleration of the instanton scheme relying heavily on the LP structure of the decoding. Similarly,  we anticipate that significant acceleration for the instanton search procedure exploring the LP structure of the load shedding formulation described above (including faster instanton-search algorithm) can be used to analyze a much larger networks, including aggregated models of the Electric Reliability Council of Texas (ERCOT) or even more ambitiously of the Eastern Grid and Western Grid of US.

\item The main point of the manuscript was to illustrate utility of the instanton-search approach in discovering weak points of the grid. To achieve this goal, while having access to only  limited observational data, we have used some very model assumptions. In particular, we have assumed Gaussian and moreover spatially uncorrelated statistics for distribution of demands, even though in reality the demands are highly correlated, in particular to external weather conditions. Our instanton-amoeba scheme is highly universal, and as such allows straightforward modification to account for effects of load correlations in future studies.

\item The optimization problem we end up solving in the search for failures is not convex and the local search algorithm we employed gives only lower bound on the probability of the most probable failure mode. It will be important to improve this aspect of the search in the future providing techniques capable to upper bound the probability and thus giving absolute guarantees. we believe that such upper bounds guarantees can be achieved with the help of the ``convexify-and-iterate" approach of the type discussed in \cite{08CWB}.

\item Sparsity of the resulting instantons became an interesting outcome of our analysis. We have seen it persistently in the tests performed, and concluded that it is robust at least with respect to the change of demand distribution. An important problem left for future studies is to confirm or disprove the robustness in broader contexts such as over other grids (graphs), other solvers (e.g. AC instead of DC), more realistic models of the demand distribution (e.g. taking into account correlations due to weather conditions and variability in fluctuations at different loads available via load forecast), other probabilistic models of fluctuations  (e.g. accounting for fluctuations in renewable generation), and other optimization techniques employed to discover the instantons.

\item Many important aspects of power flow optimization, such as these related to optimal generator dispatch differentiating between the cost of production at different power plants, were not considered in the simplified load shedding model discussed in this manuscript. To account for variability in the cost of generation should be an interesting venue for future instanton analysis.

\item Finally, this instanton-search technique can guide further exploration of more difficult statistical problems in power grid studies containing the $LP_{DC}$, or any of its inner-loop optimization substitutes, as a sub-task.  For example,  we anticipate that some version of the instanton-search technique will be useful for resolving such important problems as interdiction over power networks \cite{09SWB,10PMDL,08BV}, optimal switching of power lines to satisfy the aforementioned combinations of phase and capacity constraints \cite{09HOFO}, analysis of cascading power outages assuming sequential breakdown of saturated lines \cite{05CTD}, and also developing a reduced description of the graph to help monitor grids under stress, e.g. complementing and extending the notion of cut-sets discussed in \cite{10DPC}. We also expect that improved approaches will help to form a solid bridge to the macroscopic (and so far largely phenomenological) analysis of cascading power blackouts \cite{07DCLN,09Bar}.
\end{itemize}

\section*{Acknowledgment}

We thank all the participants of the ``Optimization and Control for Smart Grids" LDRD DR project at Los Alamos and Smart Grid Seminar Series at CNLS/LANL for multiple fruitful discussions, and especially to S. Backhaus and K. Turitsyn for very valuable comments. We are also grateful to three Referees for their very useful critique and comments and especially in helping us to connect better to realm of power engineering. Research at LANL was carried out under the auspices of the National Nuclear Security Administration of the U.S. Department of Energy at Los Alamos National Laboratory under Contract No. DE C52-06NA25396. This work was funded in part by DTRA/DOD under the grant BRCALL06-Per3-D-2-0022 on ``Network Adaptability from WMD Disruption and Cascading Failures". MC also acknowledges partial support of NMC via the NSF collaborative grant CCF-0829945  on ``Harnessing Statistical Physics for Computing and Communications''.



\end{document}